\documentclass[12pt,a4paper]{amsart}
\usepackage{amsmath,amsthm,amssymb}

\newcommand{\comment}[1]{}
\newcommand{\Floor}[1]{{\left\lfloor{#1}\right\rfloor}}

\newtheorem{theorem}{\bf Theorem}

\newtheorem{corollary}[theorem] {\bf Corollary}
\newtheorem{lemma}[theorem] {\bf Lemma}
\newtheorem{example}{\bf Example}

\newtheorem{remark}{\bf Remark}

\newcommand\PP{\mathcal P}
\newcommand\Pn{{\mathcal {P}}_n}
\newcommand\Pnc{{\mathcal {P}}^c_n}
\newcommand\LL{\mathcal L}

\newcommand{\NN}{\mathbb N}

\newcommand{\RR}{\mathbb R}
\newcommand{\CC}{\mathbb C}

\newcounter{oldresult}
\def\theoldresult{\Alph{oldresult}}
\newenvironment{oldresult}{
  \em
  \vskip 0.10in
  \refstepcounter{oldresult}
  \noindent{\bf Theorem\ \theoldresult.}
}{\vskip 0.10in}

\newcounter{oldprop}
\def\theoldprop{\Alph{oldprop}}

\newcounter{oldlemma}
\def\theoldlemma{\Alph{oldlemma}}

\newcounter{oldcor}
\def\theoldcor{\Alph{oldcor}}

\newcounter{oldconjecture}
\def\theoldconjecture{\Alph{oldconjecture}}

\newcounter{rem}
\newcounter{rev}
\usepackage{geometry,amssymb}
\begin{document}

\title[Schur type inequalities]
{Schur type inequalities for complex polynomials with no zeros in
the unit disk}

\author[Sz. Gy. R\'ev\'esz]
{Szil\'ard Gy. R\'ev\'esz}

\let\oldfootnote\thefootnote
\def\thefootnote{}

\comment{\footnotetext{Supported in part in the framework of the
Hungarian-Spanish Scientific and Technological Governmental
Cooperation, Project \# E-38/04.}}

\footnotetext{The author was supported in part by the Hungarian
National Foundation for Scientific Research, Project \#
T-049301 and K-61908.} 

\footnotetext{This work was accomplished during the author's stay
in Paris under his Marie Curie fellowship, contract \#
MEIF-CT-2005-022927.}

\let\thefootnote\oldfootnote

\address{A. R\'enyi Institute of Mathematics \newline \indent
Hungarian Academy of Sciences, \newline \indent Budapest, P.O.B.
127, 1364 Hungary.} \email{revesz@renyi.hu}

\address{and}

\address{Institut Henri Poincar\'e \newline \indent
11 rue Pierre et Marie Curie \newline \indent 75005 Paris, France}
\email{Szilard.Revesz@ihp.jussieu.fr}

\maketitle


\begin{abstract}
Starting out from a question posed by T. Erd\'elyi and J.
Szabados, we consider Schur-type inequalities for the classes of
complex algebraic polynomials having no zeroes within the unit
disk $D$.

The class of polynomials with no zeroes in $D$ -- also known as
Bernstein- or Lorentz-class -- was studied in detail earlier. For
real polynomials utilizing the Bernstein-Lorentz representation as
convex combinations of fundamental polynomials
$(1-x)^k(1+x)^{n-k}$, G. Lorentz, T. Erd\'elyi and J. Szabados
proved a number of improved versions of Schur- (and also
Bernstein- and Markov-) type inequalities.

Here we investigate the similar questions for complex polynomials.
For complex polynomials the above convex representation is not
available. Even worse, the set of complex polynomials, having no
zeroes in the unit disk, does not form a convex set. Therefore, a
possible proof must go along different lines. In fact such a
direct argument was asked for by Erd\'elyi and Szabados already
for the real case.

The sharp forms of the Bernstein- and Markov- type inequalities
are known, and the right factors are worse for complex
coefficients than for real ones. However, here it turns out that
Schur-type inequalities hold unchanged even for complex
polynomials and for all monotonic, continuous weight functions. As
a consequence, it becomes possible to deduce the corresponding
Markov inequality from the known Bernstein inequality and the new
Schur type inequality with logarithmic weight.
\end{abstract}

\bigskip
\bigskip




\section{Introduction}\label{sec:intro}

Let $\Pn$ and $\Pnc$ denote the set of one variable algebraic
polynomials of degree at most $n$ with real, resp. complex
coefficients, and denote the set of all the (real or complex)
polynomials by $\PP$ and $\PP^c$, resp.. The open unit interval
will be denoted by $I:=(-1,1)$, and the open unit disk
$\{z:\,|z|<1\}$ will be denoted by $D$. We take
\begin{equation}\label{fnormdef}
\|f\|:=\sup_{I} |f|
\end{equation}
as the norm of a polynomial or a continuous function.

In approximation theory Schur and Bernstein type polynomial
inequalities constitute an important subject, see e.g. \cite{BE,
MMR}. The classical inequality of Schur states that
\begin{equation}\label{Schurclasssic}
\|p\|\leq (n+1)\|p(x)\sqrt{1-x^2}\|\qquad\qquad \left(p\in \Pn
\right).
\end{equation}
This can be generalized to weights $(1-x^2)^{\alpha}$ with
$\alpha>0$ as well:
\begin{equation}\label{Schuralpha}
\|p\|\leq C(\alpha) n^{2\alpha}
\|p(x)(1-x^2)^{\alpha}\|\qquad\qquad \left(p\in \Pn \right).
\end{equation}
Schur's inequality \eqref{Schurclasssic} is usually combined with
Bernstein's inequality
\begin{equation}\label{Bernsteinclassic}
|p'(x)| \leq \frac{n}{\sqrt{1-x^2}} \|p\|\qquad\qquad \left(p\in
\Pn \right)
\end{equation}
to deduce Markov's inequality
\begin{equation}\label{Markovclassic}
\|p'\|\leq n^2 \|p\|\qquad\qquad \left(p\in \Pn \right).
\end{equation}

Not only Markov's inequality, but also many other results hinge
upon the basic inequalities of Schur and Bernstein. Thus there is
a well founded interest in improved versions or sharpened
inequalities of Schur and Bernstein type for various subclasses of
polynomials. An important class of interest is the Bernstein
polynomials of fixed sign, that is, the so-called ``Lorentz class"

\begin{equation}\label{Lorentzclassdef} \LL:=\{p\in\Pn ~:~
p(x)\ne 0 ~(x\in I)\}.
\end{equation}

Our interest here is the Schur type inequality for the Lorentz
class.


\section{Previous results for the Bernstein-Lorentz class}\label{sec:previres}

For $p\in\LL$, that is for real polynomials $p$ strictly positive
(or strictly negative) on the open unit interval $I:=(-1,1)$, a
so-called Lorentz representation is possible, see, e.g.,
\cite[vol. II p. 83, Aufgabe 49]{PSz}. Actually, G. Lorentz
\cite{L} considered polynomials having the representation
\begin{equation}\label{Lorentzrepi}
p(x)=\sum_{k=0}^d a_k (1-x)^k(1+x¢)^{d-k} \qquad\qquad \big(
a_k\geq 0\quad (k=1,\dots,d)\big),
\end{equation}
where $d\in \NN$ could be any natural number depending on
$p\in\PP$. Polynomials of this type were used by Lorentz \cite{L}
and others in various questions of approximation theory such as
approximation by incomplete polynomials, shape preserving
approximation and polynomials with integer coefficients. Regarding
these we refer to \cite{BE, E, K, L, KLS1, KLS2, MMR} and the
references therein.

The study of the Lorentz class \eqref{Lorentzclassdef}, the
Lorentz representation \eqref{Lorentzrepi} and the ``Lorentz
degree" $d=d(p)$  -- defined as the minimal possible degree $d$ of
such a representation of the polynomial, -- is connected to
another basic area of interest. Namely, the general idea behind
the representation \eqref{Lorentzrepi} is to exhibit the
nonnegative polynomial $p\in\Pn$ as the \emph{positive}
(nonnegative) linear combination of \emph{positive} (nonnegative)
polynomials $q_k^d(x):=(1-x)^k(1+x¢)^{d-k}$.

The positive elements $q_k^d$ form a basis of $\PP_d$, and
\eqref{Lorentzrepi} is a \emph{positive representation}, i.e., a
representation with all coefficients $a_k\geq 0$. Do all $p\geq
0$, $p\in \PP_d$ have a positive representation
\eqref{Lorentzrepi}? It is easy to see that the answer to this
question is negative. However, such questions lead to other
interesting problems, and the whole issue is a vast field of
investigations embedded into the theory of Banach lattices and
positive basis, see e.g. \cite{P1, P2, P3}. In particular, these
general results show that $\Pn$ does \emph{not} have a positive
basis at all, and, moreover, any subspace of $\Pn$ with a positive
basis has dimension at most $\Floor{n/2}$. For these questions we
refer the reader to \cite{FMPR}.

Another related matter is the theory of positive operators, in
particular, Bernstein operators
\begin{equation}\label{Bernsteinop}
B_n(f,x):=\sum_{k=0}^n f\left(\frac{2k-n}{n}\right) \binom{n}{k}
(1-x)^k(1+x¢)^{n-k}.
\end{equation}

Clearly, $B_n$ maps $C(I)$ to $\Pn$, and for $f\geq 0$ $B_n(f)\geq
0$, i.e., $B_n(f)\in \LL^{+}$, where $\LL^{+}:=\{p\in \PP~:~
p|_I>0\}$. The Bernstein operators are used extensively in the
theory of approximation, in particular for their 
shape preserving properties. 

Were now $p\in \Pn$, $p\geq 0$ a fixed point of $B_n$, comparing
\eqref{Lorentzrepi} and \eqref{Bernsteinop} would give
$p\in\LL^{+}$ and $d(p)\leq n$. Since not all $p\in \Pn\cap
\LL^{+}$ have Lorentz degree $d(p)\leq n$, we see that $B_n|_{\Pn
\cap \LL^{+}}$ can not be identity. In other words, it turns out
that the Bernstein operator is \emph{not a projection} on the set
$\Pn$. This in turn explains the shortcomings with respect to
order of approximation compared to projective operators (like,
e.g., the de la Vall\'ee Poussin operator).

Erd\'elyi and Szabados proved Schur and Bernstein type
inequalities for these polynomials using their Lorentz degree
instead of the ordinary algebraic degree. That brings into focus
the question of determining, or at least estimating the Lorentz
degree.

However, estimating the Lorentz degree of a polynomial $p\in \LL$
is usually a complicated matter. There are estimates of $d(p)$ in
terms of the zero-free region of $p$ described in \cite{ESZ} and
\cite{E}. Here we restrict our attention to the most appealing
result of this type, attributed to Lorentz, see \cite{Sc}
and \cite{ESZ}.

\begin{oldresult}{\bf (Lorentz).}\label{th:degreeq} Let $p\in\LL$.
If $p|_D\ne 0$, then we have $d(p)=\deg(p)$, the ordinary degree.
\end{oldresult}

The reason to pursue estimates of the Lorentz-degree is that there
are variants of Schur's (and also Bernstein's and Markov's)
inequalities to Lorentz polynomials with the Lorentz degree taking
over the role of the ordinary algebraic degree. Erd\'elyi and
Szabados \cite{ESZ} (see also \cite[E.14, p. 436]{BE}) have proved

\begin{oldresult}{\bf (Erd\'elyi-Szabados).}\label{th:LorentzSchur}
Let $p\in\LL$ have Lorentz degree $d=d(p)$. Then for any
$\alpha>0$ we have
\begin{equation}\label{LorentzSchur}
\|p\| \leq
\frac{(d+2\alpha)^{d+2\alpha}}{(4\alpha)^{\alpha}(d+\alpha)^{d+\alpha}}
\|p(x)(1-x^2)^{\alpha}\|\qquad \big(p\in \PP\cap
\LL,~~d=d(p)\big).
\end{equation}
\end{oldresult}

Observe that here the ``Schur constant" is of the order
$d^{\alpha}$, and in case $\alpha=1/2$ it becomes $\sqrt{d}$,
which is a considerable improvement compared to
\eqref{Schurclasssic} provided $d$ is not much larger than $n$. In
particular, combining Theorem \ref{th:degreeq} and Theorem
\ref{th:LorentzSchur} gave to Erd\'elyi and Szabados \cite{ESZ}
the following
\begin{oldresult}{\bf (Erd\'elyi-Szabados).}\label{th:specSchur}
Let $p\in \LL \cap \Pn$ and assume that $p|_D\ne 0$. Then for any
$\alpha >0$ we have
\begin{equation}\label{specSchur}
\|p\| \leq
\frac{(n+2\alpha)^{n+2\alpha}}{(4\alpha)^{\alpha}(n+\alpha)^{n+\alpha}}
\|p(x)(1-x^2)^{\alpha}\|.
\end{equation}
\end{oldresult}

Erd\'elyi and Szabados exhibit the sharpness of \eqref{specSchur}
as well. They also note that their method is bound to use
positivity of $p\in\LL$ and the result of Theorem \ref{th:degreeq}
for the Lorentz degree, while formally their end result does not
refer to Lorentz degree at all: the formulation of their results
on these inequalities does not even need the notion of Lorentz
degree and Lorentz representation for this special subclass. Hence
they comment: ``A direct proof of this statement would be
interesting."

\section{Results}\label{sec:results}

Here we will show that it is possible to obtain Theorem
\ref{th:specSchur} directly, using only nonvanishing of $p$ on
$D$. Moreover, we investigate the similar questions for complex
polynomials, where the above convex representation is not
available. It turns out that the Schur-type inequalities extend to
the complex case unchanged for all $p\in\Pnc$ (and thus without
assuming any positivity property at all), with the only assumption
of non-vanishing in $D$. This is somewhat unexpected, as an
example of Hal\'asz already established that as regards Bernstein
and Markov type inequalities, only worse estimates can be obtained
for complex polynomials \cite{E2}, \cite[p. 447]{BE}.

We formulate
\begin{theorem}\label{th:complexSchur}
Let $\varphi(t):[0,1]\to(0,\infty)$ be any decreasing, continuous
weight function. Consider a polynomial $p\in\Pnc$ and suppose that
$p|_D\ne 0$. Then
\begin{equation}\label{generalSchur}
\|p(x)\| \leq \frac{2^n}{(1+a)^n\varphi(a)} \|p(x)\varphi(x)\|=
\frac{2^n}{\|(1+x)^n\varphi(x)\|} \|p(x)\varphi(x)\|
\end{equation}
with $a\in[0,1]$ being any point of maximum of $\varphi(t)(1+t)^n$
on $[0,1]$. Moreover, equality occurs only for the polynomials
$p(x)=c(1\pm x)^n$ with $c\in\CC$ arbitrary.
\end{theorem}

\begin{corollary}\label{th:complexSchuralpha}
Let $p\in\Pnc$ and suppose that $p|_D\ne 0$. Then
\eqref{specSchur} holds true for any parameter $\alpha >0$.
Moreover, equality occurs only for the polynomials $p(x)=c(1\pm
x)^n$ with $c\in\CC$ arbitrary.
\end{corollary}


\section{Proof of Theorem \ref{th:complexSchur}}\label{sec:proof}

\begin{lemma}\label{l:factoresti}
For arbitrary $z\notin D$ and $0<a<1$ we have
\begin{equation}\label{factoresti}
\left| \frac{z-x}{z-a}\right|\leq \frac{2}{1+a} \qquad\qquad
\left( \forall x \in [a,1]\right)
\end{equation}
Moreover, equality can occur in \eqref{factoresti} only if $z=-1$
and $x=1$.
\end{lemma}
\begin{proof}
In case $\Re z \geq \frac{1+a}{2}$ we have $|z-x|\leq
\max_{x\in[a,1]}|z-x|=\max\left(|z-a|,|z-1|\right)=|z-a|$, because
for $\Re z \geq \frac{1+a}{2}$ also $|z-a|\geq |z-1|$ holds. Hence
in this case \eqref{factoresti} follows even with
$1<\frac{2}{1+a}$ on the right hand side.

In case $\Re z < \frac{1+a}{2}$ we have similarly to the above
$|z-x|\leq |z-1|$. Let us consider now the map
$f(z):=\frac{z-1}{z-a}$. This is a rational linear map of
$\widehat{\CC}\rightarrow\widehat{\CC}$ assuming real values on
$\RR$, hence is also symmetric to the real axis. Moreover, $f$
maps the set of all circles and lines to itself, $f(\infty)=1$,
$f(1)=0$, $f(a)=\infty$ and $f(-1)=\frac{2}{1+a}$. It follows that
the image of the unit circle $C=\partial D$ will be the circle $K$
symmetric to $\RR$ and going through the points $0$ and
$\frac{2}{1+a}$, that is, the circle with center $\frac{1}{1+a}$
and radius $r:=\frac{1}{1+a}$. Moreover, the domain outside of $D$
is mapped onto the interior disk $B$ of $K=\partial B$, since
$f(\infty)=1\in (0,\frac{2}{1+a}) \subset B$. However, $B\subset
D(0,\frac{2}{1+a})$, the disk centered at the origin and of radius
$\frac{2}{1+a}$. Thus for all $z\notin D$ the image satisfies
$f(z)\in B$ and therefore $|f(z)|\leq \frac{2}{1+a}$.
Consequently, we conclude in this case again that
$$
\left|\frac{z-x}{z-a}\right|\leq\left|\frac{z-1}{z-a}\right|
=|f(z)|\leq \frac{2}{1+a}.
$$
Moreover, in case $\Re z \geq \frac{1+a}{2}$ there holds a strict
inequality, and in case $\Re z < \frac{1+a}{2}$ $|z-x|=|z-1|$
entails $x=1$, and $|f(z)|=\frac{2}{1+a}$ entails $z=-1$. Thus the
assertion regarding case of equality follows, too.
\end{proof}

\begin{proof}[Proof of Theorem \ref{th:complexSchur}]
Take any parameter $0<a<1$, and consider the polynomial
\begin{equation}\label{Pdef}
P_n(x):=(1+x)^n.
\end{equation}
Plainly, for any $p(x)=\prod_{j=1}^n(x-z_j)$, where for all
$j=1,\dots,n$ we have $|z_j|\geq 1$, we have
\begin{align}\label{pestimate}
\sup_{x\in[a,1]} |p(x)| &= \sup_{a\leq x \leq 1} \left|\prod_{j=1}^n
(x-z_j)\right| = \sup_{a\leq x \leq 1}
\prod_{j=1}^n \left|\frac{x-z_j}{a-z_j}\right|\cdot |p(a)| \notag \\
&\leq \left(\frac{2}{1+a}\right)^n
\left|p(a)\right|=\frac{P_n(1)}{P_n(a)}\left|p(a)\right|~,
\end{align}
hence
\begin{equation}\label{pfinal}
\sup_{x\in[a,1]} |p(x)| \leq
\frac{P_n(1)}{P_n(a)\varphi(a)}\left|p(a)\varphi(a)\right| \leq
\frac{P_n(1)}{P_n(a)\varphi(a)}\left\|p(x)\varphi(x)\right\|~.
\end{equation}
On the other hand, for $0\leq x \leq a$ we trivially have
\begin{equation}\label{potherside}
|p(x)| \leq \frac{1}{\varphi(a)} \left|\varphi(x)p(x)\right| \leq
\frac{1}{\varphi(a)}\left\|p(x)\varphi(x)\right\|~.
\end{equation}
Combining \eqref{pfinal} and \eqref{potherside} we obtain
\begin{equation}\label{pbest}
\sup_{x\in[0,1]} |p(x)| \leq
\frac{P_n(1)}{P_n(a)\varphi(a)}\left\|p(x)\varphi(x)\right\|~,
\end{equation}
and applying this also to $p(-x)$ we finally get
\begin{equation}\label{pnormesti}
\|p\| \leq
\frac{P_n(1)}{P_n(a)\varphi(a)}\left\|p(x)\varphi(x)\right\|~.
\end{equation}
Note that \eqref{pnormesti} actually means also
\begin{equation}\label{classup}
\max\left\{ \frac{\|p\|}{\|p(x)\varphi(x)\|}~~:~~ p\in \Pnc,~~
p|_D\ne 0\right\}= \frac{\|P_n\|}{\|P_n(x)\varphi(x)\|},
\end{equation}
because \eqref{pnormesti} holds for all $0<a<1$ and hence the
maximum can be taken all over $0<a<1$.

Suppose now that we have equality in the statement of the theorem,
that is, in \eqref{generalSchur}. Since \eqref{pnormesti} was a
consequence of \eqref{pbest} and its application to $p(-x)$, case of
equality occurs only if \eqref{pbest} holds with equality either for
$p(x)$ or for $p(-x)$. Suppose, e.g., that we have equality in
\eqref{pbest} for $p(x)$, which implies equality also in
\eqref{pestimate} and \eqref{pfinal} as well. Equality in
\eqref{pestimate} in turn yields
$|p(x_0)|=\frac{P_n(1)}{P_n(a)}|p(a)|= \left(\frac{2}{1+a}\right)^n
|p(a)|$ for the maximum point $x_0\in [a,1]$ of $p$, and now the
equality part of the assertion of Lemma \ref{l:factoresti} implies
$x_0=1$ and $z_j=-1$ ($j=1,\dots,n$) for all roots of $p$. But this
proves $p(x)=c(1+x)^n$, and in case of equality for $p(-x)$, we
similarly obtain $p(x)=c(1-x)^n$. This concludes the proof.
\end{proof}

\begin{proof}[Proof of Corollary \ref{th:complexSchuralpha}]
Computing the norm on the right hand side of \eqref{classup} for
$\varphi(x)=(1-x^2)^{\alpha}$ -- that is, equivalently, taking
$a=\frac{n}{n+2\alpha}$ in \eqref{pnormesti} -- yields
\begin{equation}\label{pnormineq}
\frac{\|p\|}{\|p(x)(1-x^2)^{\alpha}\|} \leq
\frac{2^n}{\left(1+\frac{n}{n+2\alpha}\right)^n
\left(1-\frac{n^2}{(n+2\alpha)^2}\right)^{\alpha}}=
\frac{(n+2\alpha)^{n+2\alpha}}{(n+\alpha)^{n+\alpha}(4\alpha)^{\alpha}}~,
\end{equation}
which proves \eqref{specSchur}.
\end{proof}


\section{Remarks and examples}\label{sec:remex}

Comparing our proof with that of Erd\'elyi and Szabados, we can
realize that the standard approach is to make use of the convex
combination \eqref{Lorentzrepi}. Denote the set of positive
Lorentz polynomials of Lorentz degree not exceeding $d$, or
ordinary degree not exceeding $n$ by $\LL^{d}_{+}$ and by $\LL_{+
n}$, respectively. It is obvious that $\LL^{d}_{+}$, $\LL_{+ n}$
and $\LL_{+}$ are convex sets. Using convexity of $\LL^{d}_{+}$,
that is, working out proofs for the basis functions $q_{k,d}(x)$
and then adding the results, is a convenient method for real
Lorentz polynomials. However, departing real polynomials, we
necessarily need complex coefficients, and for $\PP^c\cap\LL$
similar arguments do not work. It turns out that not even the set
\begin{equation}\label{PD}
\PP_n^c(D):=\{p\in\Pn^c ~~:~~ p|_I>0,~p|_D\ne 0 \}
\end{equation}
is convex; hence convex combinations can not be used directly in
this setting.

\begin{example} Let $0<a<1$, $w:=a+i\sqrt{1-a^2}$, and consider
the polynomials
\begin{align}\label{nonconvexex}
p(x):&=(1-x)(x^2-2ax+1)=(1-x)(x-w)(x-\overline{w})\notag \\
& =-x^3+(1+2a)x^2-(2a+1)x+1 \\
q(x):&=1+x+x^2+x^3=(1+x)(1+x^2)=(x+1)(x+i)(x-i)~ \notag.
\end{align}
Then $p,q\in \PP_3^c(D)$, but for $ r:=\frac{1}{2}p+ \frac{1}{2}q
$ one has $r\notin \PP_3^c(D)$, hence $\PP_3^c(D)$ is not convex.
\end{example}

Indeed, both $p$ and $q$ have zeroes of absolute value 1 only, so
they belong to $\PP_3^c(D)$. Moreover, for
$$
r(x)=\frac{p(x)+q(x)}{2}=(1+a)x^2-ax+1
$$
we obviously have $r\in \LL_{+}$ ($\LL_{+}$ is convex). On the
other hand the roots of $r(x)$ are
\begin{equation}\label{rrots}
x_{1,2}=\frac{a\pm\sqrt{a^2-4(1+a)}}{2(1+a)}=\frac{a\pm i
\sqrt{4+4a-a^2}}{2(1+a)}~.
\end{equation}
Observe that $4+4a-a^2> 4 >0$ for all $a\in (0,1)$. Now we can
compute
\begin{equation}\label{rootabsval}
|x_{1,2}|^2=\frac{a^2+4+4a-a^2}{4(1+a)^2}=\frac{1}{1+a}<1~,
\end{equation}
hence $|x_{1,2}|=1/\sqrt{1+a}<1$, $x_{1,2}\in D$ and $r\notin
\PP_n^c(D)$ for any $n\in\NN$.

Note that in this example both $p$ and $q$ have degree 3, and by
Theorem \ref{th:degreeq} and $p,q\in \PP_3^c(D)$ their Lorentz
degree is 3. Consequently, by convexity of $\LL_{+}$ and
$\LL^{d}_{+}$, we must have $d(r)\leq 3$, while $d(r)\geq {\rm
deg}~r =2$. To decide the exact value of $d(r)$, note that
$(1+x)^2$, $1-x^2$ and $(1-x)^2$ form a basis of $\PP_2$, and easy
linear algebra gives $r(x)=\frac12 (1+x)^2 + \frac{1+a}{2} (1-x)^2
-\frac{a}{2} (1-x^2)$, whence the unique degree 2 representation
is not positive and the Lorentz degree can not be $2$. Actually,
$d(r)\ne {\rm deg} ~r$ already follows from \cite[Theorem 2
(ii)]{ESZ} or \cite[Proposition, p. 117]{ESZ}. Whence $d(r)=3$,
and the corresponding representation is easyly obtained from those
of $p$ and $q$.

The following comment was offered by Tam\'as Erd\'elyi.

\begin{remark}[\bf Erd\'elyi]
As regards Schur's inequality, we have a better than general bound
\eqref{specSchur} at least for the class $\PP_n^c(D)$. In fact,
this can also be obtained from the real case, i.e., from Theorem
\ref{th:LorentzSchur} and \ref{th:degreeq}, independently of
Theorem \ref{th:complexSchur} or Corollary
\ref{th:complexSchuralpha}.
\end{remark}
Indeed, let $p\in\Pnc$ such that $p(z)\ne 0$ for $z\in D$.
Consider also
$$
\widetilde{p}(z):=\overline{p(\overline{z})}=\prod_{j=1}^n
(z-\overline{z_j}) \qquad\qquad \left(p(z):=\prod_{j=1}^n (z-z_j)
\right)~
$$
and take $p^{*}(z):=p(z)\widetilde{p}(z)$. Obviously $p^{*}\in
\PP_{2n}$ and $p^{*}\in \LL_{+2n}$, too. Applying Theorem
\ref{th:LorentzSchur} with power $\alpha^{*}:=2\alpha$ to $p^{*}$
of degree $n^{*}:=2n$ we get
\begin{gather*}
\|p\|^2=\|p^{*}\|=\|p\widetilde{p}\|\leq
\frac{(2n+4\alpha)^{2n+4\alpha}}{(8\alpha)^{2\alpha}(2n+2\alpha)^{2n+2\alpha)}}
\|p(x)\widetilde{p}(x)(1-x^2)^{2\alpha}\| \\ = \left(
\frac{(n+2\alpha)^{n+2\alpha}}{(4\alpha)^{\alpha}(n+\alpha)^{n+\alpha}}
\|p(x)(1-x^2)^{\alpha}\| \right)^2~.
\end{gather*}

However, for the Bernstein and Markov inequalities in the
generality of complex polynomials not vanishing in $D$, we have
substantially worse factors, see \cite[p.474]{BE} and \cite{E2}.
The example of Hal\'asz below shows what we can expect at most.

\begin{example}[\bf Hal\'asz] Let $m$ be chosen as $[(n-1)/2]$, so that
$2m+1\leq n \leq 2m+2$. Define the $\deg P=n$ polynomial $P$ as
$$
P(z):=(z-1)\prod_{j=1}^m \left(z-e^{\frac{2\pi
ij}{2m+1}}\right)^2.
$$
Then $P|_D\ne 0$, $\|P\|_D=2=|P(-1)|$ and $P'(-1) \gg c n \log n$.
Moreover, for any $x\in[-1,1]$, we also have $P'(x)> c n \log
\frac{e}{1-x^2}$ whenever this is smaller than $cn\log n$.

Consequently, no better bound, than $c \min \left( n\log n; n \log
\frac{e}{1-x^2}\right)$ is valid in the Markov- and Bernstein
inequality, even if we restrict to $\PP_n^c(D)$.
\end{example}

The (essentially standard) calculation showing these lower
estimates can be found, e.g., in \cite{E2} or \cite[p. 447]{BE}.
These are indeed the right factors as the corresponding upper
estimation is proved, e.g., in \cite{E2}.

A standard way of proving Markov type inequalities is to combine
Bernstein inequalities with Schur inequalities. Of course, to get
a sharp Markov estimate we must combine sharp Bernstein and sharp
Schur inequalities as well. Thanks to the general form (with any
monotone $\varphi(x)$) of our formulation of the Schur type
inequality Theorem \ref{th:complexSchur}, here we can indeed
deduce the Markov bound from the corresponding Bernstein
inequality. Indeed, the known Bernstein type estimate (see
\cite[Theorem 2.1]{E2}) says
\begin{equation}\label{complexbernstein}
|p'(x)| \leq n\log\frac{e}{1-x^2}\qquad (|x|<1,~~ p\in
\PP_n^c(D)),
\end{equation}
and applying the Schur inequality \eqref{generalSchur} to $p'(x)$
and $\varphi(x):= \log^{-1}\frac{e}{1-x^2}$ we obtain
$$
\|p'(x)\| \leq \frac{2^n}{\|(1+x)^n \varphi(x)\|} \|\varphi(x)
p'(x) \| \leq \frac{2^n}{|(1+x_0)^n \varphi(x_0)|} \|\varphi(x)
p'(x) \|
$$
with arbitrary $x_0\in I$. Choosing $x_0:=1-2/n$, say, we thus
obtain $ \|p'\| \leq C \log n \| \varphi p'\|$ and this can be
estimated by the above Bernstein inequality
\eqref{complexbernstein} as $\leq C n \log n $.

Note that given the logarithmic weight in the complex case,
restricting to weights $(1-x^2)^{\alpha}$ would bring by itself
the loss of the possibility of this deduction.


\end{document}